\documentclass[12pt,twoside]{article} 

\usepackage{amssymb,amsmath}

\date{} 

\begin{document} 

\centerline{\bf Applied Mathematical Sciences, Vol. x, 20xx, no. xx, xxx - xxx} 

\centerline{\bf HIKARI Ltd, \ www.m-hikari.com} 

\centerline{\bf http://dx.doi.org/xx.xxxxx/}

\centerline{} 

\centerline{} 

\centerline{\Large{\bf Global Well-Posedness}} 

\centerline{} 

\centerline{\Large{\bf for a 2-D Viscoelastic Fluid Model}} 


\centerline{} 

\centerline{} 

\centerline{\bf {Mikhail A. Artemov}} 

\centerline{} 

\centerline{Voronezh State University} 

\centerline{Universitetskaya Pl. 1, Voronezh, Russian Federation, 394006} 

\centerline{E-mail address: artemov\_m\_a@mail.ru} 
\centerline{} 

\centerline{\bf {George G. Berdzenishvili}} 

\centerline{} 

\centerline{Voronezh State University} 

\centerline{Universitetskaya Pl. 1, Voronezh, Russian Federation, 394006} 

\centerline{E-mail address: berdzen\_g@mail.ru} 

\newtheorem{Theorem}{\quad Theorem}

\newtheorem{Definition}{\quad Definition} 

\newtheorem{Corollary}[Theorem]{\quad Corollary} 

\newtheorem{Lemma}{\quad Lemma} 

\newtheorem{Example}[Theorem]{\quad Example} 

\centerline{}

{\footnotesize Copyright $\copyright$ 20xx Mikhail A. Artemov. This article is distributed under the Creative Commons Attribution License, which permits unrestricted use, distribution, and reproduction in any medium, provided the original work is properly cited.}

\begin{abstract} This paper is concerned with a mathematical model which describes 2-D flows of an incompressible viscoelastic fluid of Oldroyd type in a bounded domain. We prove the existence and uniqueness theorem for global (in time) weak solutions and derive the energy equation.

\end{abstract} 



\section{Introduction} In this paper, we will consider the following nonlinear system, which describes, at least at first approximation, 2-D flows of an incompressible viscoelastic fluid of the Oldroyd kind in a fixed geometry $\Omega$:
\begin{equation}\label{equ-1}
\text{Re}\left(\dfrac{\partial \boldsymbol{v}}{\partial t}+(\boldsymbol{v}\cdot\nabla)\boldsymbol{v}\right)-(1-a)\Delta\boldsymbol{v}-\nabla\cdot\boldsymbol{\tau}+\nabla \pi=\boldsymbol{f}\;\; \text{in}\;\;\Omega\times(0,T), 
\end{equation}
\begin{equation}\label{equ-2}
\nabla\cdot\boldsymbol{v}=0\;\; \text{in}\;\; \Omega\times(0,T),
\end{equation}
\begin{equation}\label{equ-3}
\text{We}\dfrac{\partial\boldsymbol{\tau}}{\partial t}+\boldsymbol{\tau}=2a\boldsymbol{\mathcal E}(\boldsymbol{v})\;\;\text{in}\;\; \Omega\times(0,T),
\end{equation}
\begin{equation}\label{equ-4}
\boldsymbol{v}=\boldsymbol{0}\;\;\text{on}\; \partial\Omega\times(0,T), 
\end{equation}
\begin{equation}\label{equ-6}
\boldsymbol{v}(\cdot, 0)=\boldsymbol{v}_0, \;\;\boldsymbol{\tau}(\cdot, 0)=\boldsymbol{\tau}_0  \;\;\text{in}\;\Omega.
\end{equation}
Here, $\Omega$ is an open bounded (sufficiently smooth) subset of $\mathbb{R}^2$,  $T$ is a prescribed final time, $\boldsymbol{v}=\{v_{i}(\boldsymbol{x},t)\}$ is the velocity of the fluid, $\boldsymbol{\tau}=\{\tau_{ij}(\boldsymbol{x},t)\}$ is the elastic extra-stress tensor with $\tau_{12}(\boldsymbol{x},t)\equiv\tau_{21}(\boldsymbol{x},t)$, $\pi=\pi(\boldsymbol{x},t)$ is the pressure distribution, $\boldsymbol{\mathcal E}(\boldsymbol{v})$ is the symmetric part of the velocity gradient, 
$$
\boldsymbol{\mathcal E}(\boldsymbol{v})=\frac{1}{2}\left(\nabla \boldsymbol{v}+(\nabla \boldsymbol{v})^{\rm T}\right),
$$
 the vector function $\boldsymbol{f}=\{f_i(\boldsymbol{x},t)\}$ is the field of exterior forces,  $\text{Re}>0$ is the Reynolds number of the fluid,  $\text{We}>0$ is the Weissenberg number, and $a$ is a model parameter such that $0<a<1$.   

In system \eqref{equ-1}--\eqref{equ-6}, the unknowns are $\boldsymbol{v}$, $\boldsymbol{\tau}$ and $p$; the data of the problem are $\boldsymbol{f}$, $\boldsymbol{v}_0$, $\boldsymbol{\tau}_0$ and the constants  $\text{Re}$, $\text{We}$ and $a$. 

 It should be mentioned at this point that a lot of mathematical studies have been conducted towards the Oldroyd systems and similar non-Newtonian models.  Starting with Renardy \cite{ren85} and Guillop\'e and Saut \cite{guisau90} there is an ever growing list of contributions. Local existence and uniqueness theorems are given in \cite{bes01}, \cite{fergilort98}, \cite{guisau90}. In his paper \cite{tur95}, Turganbaev proves existence of global (in time) weak solutions for Oldroyd fluids with material time derivative in the constitutive law. Various well-posedness results for motion equations with classical boundary conditions were obtained also in \cite{agr98}, \cite{bar13}, \cite{chemia08}, \cite{liomas00}. Note also that slip problems for Oldroyd type models have recently attracted considerable interest (see \cite{art-e}, \cite{art-a}, \cite{bar14-s}, \cite{bar14-a}, \cite{rou14}). A review of mathematical results on viscoelastic fluids models and many references may be found in \cite{sau13}.
  
 Unfortunately, there exists only a limited number of results concerning with the uniqueness of solutions for Oldroyd type models, and all of them are devoted to sufficiently slow flows or local-in-time solutions (see \cite{bes01}, \cite{chemia08}, \cite{guisau90}, \cite{fergilort98}). A possible way out is to consider 2-D flows with linearized constitutive law~\eqref{equ-3}. So, we arrive at problem \eqref{equ-1}--\eqref{equ-6}.

Our main result provides the existence and uniqueness of weak solutions (see Theorem~\ref{mainthr1}). 
The proof of this result is based on classical (in mathematical studies of hydrodynamics) techniques, including the Faedo--Galerkin procedure, energy estimates and Ladyzhenskaya's inequality.

\section{Preliminary Notes} 
Throughout this paper we will use the following spaces: the Lebesgue spaces $\boldsymbol{L}_p(\Omega)$, $1\leq p\leq\infty$, with the norm $\|\cdot\|_{\boldsymbol{L}_p}$; the Sobolev spaces $\boldsymbol{W}^m_q(\Omega)$, $m\in\mathbb{N}$, $1\leq q\leq\infty$, with the norm $\|\cdot\|_{\boldsymbol{W}^m_q}$.  This notation is used for vector or matrix functions, and which of them is clear from the context. As usual,  
$$
\boldsymbol{H}^m(\Omega):=\boldsymbol{W}^m_2(\Omega).
$$
See \cite{ada03} for definitions and properties of these spaces.
 
We will use  the standard (in mathematical studies of hydrodynamics) function spaces: 
\begin{align*}
&\mathcal{D}(\Omega):=\{\boldsymbol{u}\in \boldsymbol{C}^{\infty}(\Omega):\;\mathop{\rm supp}\,\boldsymbol{u}\subset\Omega\},
\\
&\mathcal{V}:=\{\boldsymbol{u}\in \boldsymbol{C}^{\infty}(\Omega):\;\nabla\cdot\boldsymbol{u}=0,\;  \mathop{\rm supp}\,\boldsymbol{u}\subset\Omega\},
\\
&\boldsymbol{H}:=\text{the closure of}\;\;\mathcal{V}\;\; \text{with respect to the norm of }\; \boldsymbol{L}_{2}(\Omega),
\\
&\boldsymbol{V}:=\text{the closure of}\;\;\mathcal{V}\;\; \text{with respect to the norm of }\; \boldsymbol{H}^{1}(\Omega),
\\
&\boldsymbol{H}^m_0(\Omega):=\text{the closure of}\;\;\mathcal{D}(\Omega)\;\; \text{with respect to the norm of }\; \boldsymbol{H}^{m}(\Omega).
\end{align*}

We introduce the following notation
$$
(\boldsymbol{u},\boldsymbol{w}):=\int\limits_\Omega\boldsymbol{u}\cdot\boldsymbol{w}\,dx,\;\;\;\boldsymbol{u},\boldsymbol{w}\in \boldsymbol{L}_{2}(\Omega).
$$

Let us define the inner product in $\boldsymbol{V}$ by the formula
$$
(\boldsymbol{u},\boldsymbol{w})_{\boldsymbol{V}}:=(1-a)(\nabla\boldsymbol{u},\nabla\boldsymbol{w}).
$$
This inner product is convenient for deriving estimates of solutions.
The norm $\|\boldsymbol{u}\|_{\boldsymbol{V}}:=(\boldsymbol{u},\boldsymbol{u})^{1/2}_{\boldsymbol{V}}$ is equivalent to the norm induced from the space $\boldsymbol{H}^1(\Omega)$.  

Using the Riesz representation theorem, one identifies the space $\boldsymbol{H}$ with its conjugate space $\boldsymbol{H}^*\equiv\boldsymbol{H}$. Therefore we arrive at the inclusions
$$
\boldsymbol{V}\subset\boldsymbol{H}\equiv \boldsymbol{H}^*\subset\boldsymbol{V}^*.
$$

For a Banach space $\boldsymbol{X}$, denote by $\boldsymbol{C}([0,T];\boldsymbol{X})$ the space of continuous functions $\boldsymbol{u}\colon[0,T]\to\boldsymbol{X}$. Furthermore, denote by $\boldsymbol{C}_w([0,T];\boldsymbol{X})$ the set of the functions $\boldsymbol{u}\colon[0,T]\to\boldsymbol{X}$ which are weakly continuous.
As usual, $\boldsymbol{L}_q(0,T;\boldsymbol{X})$ denotes the space of $L_q$-integrable functions from $[0, T]$ into $\boldsymbol{X}$.

 Finally, the symbol $\mathbb{R}^{2\times 2}_{\text{sym}}$ will denote the space of symmetric real $2\times 2$-matrices.

\section{Weak Formulation of  Problem \eqref{equ-1}--\eqref{equ-6}} 

We now give the following definition. Let
\begin{equation}\label{treb}
\boldsymbol{v}_0\in \boldsymbol{H},\;\boldsymbol{\tau}_0\in \boldsymbol{L}_2(\Omega),\;\boldsymbol{f}\in \boldsymbol{L}_2\bigl(0,T;\boldsymbol{L}_2(\Omega)\bigr).
\end{equation}

\begin{Definition} 
 One says that a pair of functions
\begin{align*}
\boldsymbol{v}&\colon\Omega\times[0,T]\to\mathbb{R}^2,\;\;\;\;\;\;\boldsymbol{v}\in \boldsymbol{L}_2(0,T; \boldsymbol{V})\cap\boldsymbol{C}_w([0,T]; \boldsymbol{H}),
\\
\boldsymbol{\tau}&\colon\Omega\times[0,T]\to\mathbb{R}^{2\times 2}_{\text{sym}},\;\;\; \boldsymbol{\tau}\in \boldsymbol{C}_w([0,T]; \boldsymbol{L}_2(\Omega))
\end{align*}
is a {\it weak solution} to problem {\rm \eqref{equ-1}--\eqref{equ-6}} if it satisfies the initial conditions \eqref{equ-6}, and if the equalities
\begin{gather*}
{\rm Re}\frac{d}{dt}(\boldsymbol{v},\boldsymbol{\varphi})-{\rm Re}\sum\limits_{i=1}^{2}\Bigl(v_i\boldsymbol{v},\frac{\partial{\boldsymbol{\varphi}}}{\partial x_i}\Bigr)
+(1-a)\bigl(\nabla \boldsymbol{v},\nabla \boldsymbol{\varphi}\bigr)
\\
+\bigl(\boldsymbol{\tau},\boldsymbol{\mathcal E}(\boldsymbol{\varphi})\bigr)
=(\boldsymbol{f},\boldsymbol{\varphi})\;\;\;\forall\boldsymbol{\varphi}\in \boldsymbol{V},
\end{gather*}
$$
{\rm We}\frac{d}{dt}(\boldsymbol{\tau},\boldsymbol{\psi})
+(\boldsymbol{\tau},\boldsymbol{\psi})
=2a\bigl(\boldsymbol{\mathcal E}(\boldsymbol{v}),\boldsymbol{\psi}\bigr)\;\;\;\;\forall\boldsymbol{\psi}\in \boldsymbol{L}_2(\Omega)
$$
are true in the distribution sense on $(0, T)$.
\end{Definition} 

\section{Main Result} 

\begin{Theorem}\label{mainthr1}
For given $\boldsymbol{f}$, $\boldsymbol{v}_0$  and $\boldsymbol{\tau}_0$ which satisfy \eqref{treb}, there exists a unique weak solutions $(\boldsymbol{v},\boldsymbol{\tau})$ to problem  {\rm \eqref{equ-1}--\eqref{equ-6}}. Moreover, $\boldsymbol{v}$ belongs to $\boldsymbol{C}([0,T]; \boldsymbol{H})$, $\boldsymbol{\tau}$ belongs to $\boldsymbol{C}([0,T]; \boldsymbol{L}_2(\Omega))$, and the following energy equation holds:
\begin{gather}
{\rm Re}\|\boldsymbol{v}(t)\|_{\boldsymbol{L}_2}^2+2(1-a)\int\limits_0^t\|\nabla\boldsymbol{v}(s)\|_{\boldsymbol{L}_2}^2\,ds 
+\frac{1}{a}\int\limits_0^t\|\boldsymbol{\tau}(s)\|_{\boldsymbol{L}_2}^2\,ds\nonumber
\\
+\frac{{\rm We}}{2a}\|\boldsymbol{\tau}(t)\|_{\boldsymbol{L}_2}^2
=2\int\limits_0^t(\boldsymbol{f}(s),\boldsymbol{v}(s))\,ds+{\rm Re}\|\boldsymbol{v}_0\|_{\boldsymbol{L}_2}^2+\frac{{\rm We}}{2a}\|\boldsymbol{\tau}_0\|_{\boldsymbol{L}_2}^2.\label{eneq}
\end{gather}
\end{Theorem} 

\section{Proof of Theorem~\ref{mainthr1}} 
For reader's convenience, let us recall some results will be used later.

\begin{Lemma}[the Ladyzhenskaya inequality]\label{lm-1}
For any open set $\Omega\subset\mathbb{R}^2$ and any $w\in{H}^1_0(\Omega)$, 
$$
\|w\|_{L_4}\leq 2^{1/4}\|w\|_{L_2}^{1/2}\|\nabla w\|_{\boldsymbol{L}_2}^{1/2}.
$$
\end{Lemma}

\begin{Lemma}\label{lm-2}
Let $\boldsymbol{X}$ and $\boldsymbol{Y}$ be Hilbert spaces such that
$$
\boldsymbol{X}\subset\boldsymbol{Y}\equiv \boldsymbol{Y}^*\subset\boldsymbol{X}^*.
$$
If a function $\boldsymbol{w}$ belongs to the space $\boldsymbol{L}_2(0, T; \boldsymbol{X})$ and its derivative $\boldsymbol{w}^\prime$ belongs to $\boldsymbol{L}_2(0, T; \boldsymbol{X}^*)$, then $\boldsymbol{w}$ is almost everywhere equal to a function continuous from $[0, T]$ into $\boldsymbol{Y}$ and the following equality holds in the scalar distribution sense on $(0, T)$:
$$
\frac{d}{dt}\|\boldsymbol{w}(t)\|^2_{\boldsymbol{Y}}=2\langle\boldsymbol{w}^\prime(t),\boldsymbol{w}(t)\rangle_{\boldsymbol{X}^*\times\boldsymbol{X}}.
$$
\end{Lemma}

The proofs of Lemma \ref{lm-1} and Lemma \ref{lm-2} are given in the monograph \cite{tem77}.
\vspace{2mm}

{ \it Proof of Theorem~\ref{mainthr1}. }Let us start by proving the existence of weak solutions to problem \eqref{equ-1}--\eqref{equ-6}. In order to obtain a weak solution, we apply the Faedo--Galerkin procedure.

Let $\{\boldsymbol{\varphi}^j\}_{j=1}^\infty$ be an orthonormal basis of the space $\boldsymbol{H}$ such that the system $\{\boldsymbol{\varphi}^j\}_{j=1}^\infty$ is total in the space $\boldsymbol{V}$, and  let $\{\boldsymbol{\psi}^j\}_{j=1}^\infty$ be an orthonormal basis of the space $\boldsymbol{L}_2\bigl(\Omega)$.  

For each $n\in\mathbb{N}$ we define an approximate solution as follows:
$$\boldsymbol{v}^n=\sum_{j=1}^na_{nj}(t)\boldsymbol{\varphi}^j(\boldsymbol{x}),\;\;\boldsymbol{\tau}^n=\sum_{j=1}^nb_{nj}(t)\boldsymbol{\psi}^j(\boldsymbol{x}),$$
where $a_{nj}$ and $b_{nj}$ are unknown functions, and
\begin{gather}
\mathop{\rm Re}\Bigl(\frac{\partial \boldsymbol{v}^n}{\partial t}, \boldsymbol{\varphi}^j \Bigr)+\mathop{\rm Re}\sum\limits_{i=1}^{2}\Bigl(v_i^n\frac{\partial \boldsymbol{v}^n}{\partial x_i},\boldsymbol{\varphi}^j\Bigr)+(\boldsymbol{\tau}^n, \boldsymbol{\mathcal E}(\boldsymbol{\varphi}^j))
\nonumber
\\
+(1-a)(\nabla\boldsymbol{v}^n, \nabla\boldsymbol{\varphi}^j)
 =\left(\boldsymbol{f}, \boldsymbol{\varphi}^j \right), \;j=1,2,\dots,n,
 \label{galer-sys-1}
\end{gather}
\begin{equation}\label{galer-sys-2}
\mathop{\rm We}\Bigl(\frac{\partial {\boldsymbol{\tau}}^n}{\partial t},\boldsymbol{\psi}^j\Bigr) +(\boldsymbol{\tau}^n,\boldsymbol{\psi}^j)
=2a(\boldsymbol{\mathcal E}(\boldsymbol{v}^n),\boldsymbol{\psi}^j),  \;j=1, 2,\dots,n,
\end{equation}
\begin{equation}\label{galer-sys-4}
\boldsymbol{v}^n(\cdot,0)=\sum\limits_{j=1}^{n}(\boldsymbol{v}_0, \boldsymbol{\varphi}^j)\boldsymbol{\varphi}^j,\;\;{\boldsymbol{\tau}}^n(\cdot,0)=\sum\limits_{j=1}^{n}({\boldsymbol{\tau}}_0, \boldsymbol{\psi}^j)\boldsymbol{\psi}^j.
\end{equation}

Suppose that a pair $(\boldsymbol{v}^n,{\boldsymbol{\tau}}^n)$ satisfies \eqref{galer-sys-1}--\eqref{galer-sys-4}. We claim that for $(\boldsymbol{v}^n,{\boldsymbol{\tau}}^n)$ there exist suitable a priori estimates, which are independent of $n$. Indeed, multiplying \eqref{galer-sys-1} by $a_{nj}(t)$ and summing over $j$ from 1 to $n$, we have
\begin{equation}\label{galer-sys-1-1}
\mathop{\rm Re}\Bigl(\frac{\partial \boldsymbol{v}^n}{\partial t}, \boldsymbol{v}^n \Bigr)+\left(\boldsymbol{\tau}^n, \boldsymbol{\mathcal E}(\boldsymbol{v}^n) \right)+(1-a)\bigl(\nabla\boldsymbol{v}^n, \nabla\boldsymbol{v}^n\bigr)
=\left(\boldsymbol{f}, \boldsymbol{v}^n \right).
\end{equation}
Further, we multiply \eqref{galer-sys-2} by  $b_{nj}(t)$ add the results for $j=1,2\dots, n$:
\begin{equation}\label{galer-sys-2-2}
\mathop{\rm We}\Bigl(\frac{\partial \boldsymbol{\tau}^n}{\partial t},\boldsymbol{\tau}^n\Bigr)+ (\boldsymbol{\tau}^n,\boldsymbol{\tau}^n)=2a(\boldsymbol{\mathcal E}(\boldsymbol{v}^n),\boldsymbol{\tau}^n).  
\end{equation}
Finally, we multiply \eqref{galer-sys-1-1} by $2a$. Summing the result and $\eqref{galer-sys-2-2}$, we can infer
$$
2a\mathop{\rm Re}\Bigl(\frac{\partial \boldsymbol{v}^n}{\partial t}, \boldsymbol{v}^n  \Bigr)+2a(1-a)(\nabla\boldsymbol{v}^n,\nabla\boldsymbol{v}^n)
$$
\begin{equation}\label{est-10}
+(\boldsymbol{\tau}^n,\boldsymbol{\tau}^n)+\mathop{\rm We}\Bigl(\frac{\partial \boldsymbol{\tau}^n}{\partial t},\boldsymbol{\tau}^n\Bigr)=2a\left(\boldsymbol{f}, \boldsymbol{v}^n \right).
\end{equation} 
Using the Gr\"onwall--Bellman inequality, from \eqref{est-10} we deduce that the sequence $\{\boldsymbol{v}^n\}_{n=1}^\infty$ is bounded in the spaces $\boldsymbol{L}_2(0,T; \boldsymbol{V})$ and $\boldsymbol{L}_\infty(0,T; \boldsymbol{H})$.  Also, the sequence $\{\boldsymbol{\tau}^n\}_{n=1}^\infty$ is bounded in the space $\boldsymbol{L}_\infty(0,T; \boldsymbol{L}_2(\Omega))$. By virtue of these properties, problem \eqref{galer-sys-1}--\eqref{galer-sys-4} admits a global (in time) solution for any $n\in\mathbb{N}$. Moreover, by analogy with the case of the Navier--Stokes system (see \cite{tem77}), we can apply the compactness theorems for passing to the limit as $n\to\infty$. As a result, we get a weak solution of problem \eqref{equ-1}--\eqref{equ-6}.

We will prove  that this solution is unique. The proof is by reductio ad absurdum. Let us assume that $(\boldsymbol{v}_1,\boldsymbol{\tau}_1)$ and $(\boldsymbol{v}_2,\boldsymbol{\tau}_2)$ are weak solutions to problem  \eqref{equ-1}--\eqref{equ-6}.
We will show that $\boldsymbol{v}_1=\boldsymbol{v}_2$  and $\boldsymbol{\tau}_1=\boldsymbol{\tau}_2.$

First let us show that
\begin{equation}\label{en-1}
\boldsymbol{v}_i^\prime\in \boldsymbol{L}_2(0,T;\boldsymbol{V}^*), \;\;i=1,2. 
\end{equation}
To do so, we introduce auxiliary operators by the following formulas:
$$
\mathcal{A}\colon\boldsymbol{V}\to\boldsymbol{V}^*,\;\;
\langle\mathcal{A}(\boldsymbol{v}),\boldsymbol{\varphi}\rangle=(a-1)\bigl(\nabla\boldsymbol{v}, \nabla\boldsymbol{\varphi}\bigr),
$$
$$
\mathcal{B}\colon\boldsymbol{V}\times\boldsymbol{V}\to\boldsymbol{V}^*,\;\;
\langle\mathcal{B}(\boldsymbol{v},\boldsymbol{w}),\boldsymbol{\varphi}\rangle=\text{Re}\sum\limits_{i=1}^{2}\Bigl(v_i\boldsymbol{w},\frac{\partial\boldsymbol{\varphi}}{\partial x_i}\Bigr),
$$
$$
\mathcal{G}\colon\boldsymbol{L}_2(\Omega)\to\boldsymbol{V}^*,\;\; \langle\mathcal{G}(\boldsymbol{\tau}),\boldsymbol{\varphi}\rangle=-\bigl(\boldsymbol{\tau}, \boldsymbol{\mathcal E}(\boldsymbol{\varphi})\bigr).
$$
One can deduce from the definition of weak solutions that
$$
{\rm Re}\,\boldsymbol{v}_i^\prime=\mathcal{A}(\boldsymbol{v}_i)+\mathcal{B}(\boldsymbol{v}_i,\boldsymbol{v}_i)+\mathcal{G}(\boldsymbol{\tau}_i)+\boldsymbol{f}.
$$
It is readily seen that
\begin{equation}\label{en-2}
\mathcal{A}(\boldsymbol{v}_i),  \mathcal{G}(\boldsymbol{\tau}_i),\boldsymbol{f}\in \boldsymbol{L}_2(0,T;\boldsymbol{V}^*).
\end{equation}
We claim that
\begin{equation}\label{en-3}
\mathcal{B}(\boldsymbol{v}_i)\in \boldsymbol{L}_2(0,T;\boldsymbol{V}^*).
\end{equation}
Indeed, by the H\"older inequality and the Ladyzhenskaya inequality (Lemma~\ref{lm-1}), we obtain
\begin{align}
|\langle\mathcal{B}(\boldsymbol{v}_i(t),\boldsymbol{v}_i(t)),\boldsymbol{\varphi}\rangle|
&\leq C_1 \|\boldsymbol{v}_i(t)\|_{\boldsymbol{L}_4}^2\|\boldsymbol{\varphi}\|_{\boldsymbol{V}}\nonumber
\\
&\leq C_2\|\boldsymbol{v}_i(t)\|_{\boldsymbol{L}_2}\|\nabla\boldsymbol{v}_i(t)\|_{\boldsymbol{L}_2}\|\boldsymbol{\varphi}\|_{\boldsymbol{V}}. \label{dp-1}
\end{align}
For the rest of this article, by $C_i$, $i=1,2,\dots$, we will denote positive constants, which are independent of $t$. 
From estimate \eqref{dp-1} it follows that \eqref{en-3} holds. 

If we combine \eqref{en-3} with \eqref{en-2}, we get \eqref{en-1}.

Moving on, let us prove that
\begin{equation}\label{en-4}
\boldsymbol{\tau}_i^\prime\in \boldsymbol{L}_2\bigl(0,T;[\boldsymbol{H}^1(\Omega)]^*\bigr),\;i=1,2.
\end{equation}
Introducing the auxiliary operators
$$
\mathcal{K}\colon\boldsymbol{L}_2(\Omega)\to[\boldsymbol{H}^1(\Omega)]^*,\;\;
\langle\mathcal{K}(\boldsymbol{\tau}),\boldsymbol{\psi}\rangle=-(\boldsymbol{\tau},\boldsymbol{\psi}),
$$
$$
\mathcal{N}\colon\boldsymbol{V}\to[\boldsymbol{H}^1(\Omega)]^*,\;\;
\langle\mathcal{N}(\boldsymbol{v}),\boldsymbol{\psi}\rangle=2a\bigl(\boldsymbol{\mathcal E}(\boldsymbol{v}),\boldsymbol{\psi}\bigr),
$$
we see that
$$
{\rm We}\,\boldsymbol{\tau}_i^\prime=\mathcal{K}(\boldsymbol{\tau}_i)+\mathcal{N}(\boldsymbol{v}_i).
$$
It can easily be checked that
\begin{equation*}
\mathcal{K}(\boldsymbol{\tau}_i),\mathcal{N}(\boldsymbol{v}_i)\in \boldsymbol{L}_2(0,T;[\boldsymbol{H}^1(\Omega)]^*),
\end{equation*}
which proves that \eqref{en-4} is valid.

Using \eqref{en-1} and \eqref{en-4}, we apply Lemma \ref{lm-2} to $\boldsymbol{v}_i$ and $\boldsymbol{\tau}_i$, and we obtain 
$$
\boldsymbol{v}_i\in \boldsymbol{C}([0,T]; \boldsymbol{H}(\Omega)),\;\;\boldsymbol{\tau}_i\in\boldsymbol{C}([0,T]; \boldsymbol{L}_2(\Omega)),\;\;i=1,2.
$$
Moreover, the following equalities hold almost everywhere on $[0,T]$:
\begin{equation}\label{en-7}
\frac{d}{dt}\|\boldsymbol{u}(t)\|^2_{\boldsymbol{L}_2}=2\langle\boldsymbol{u}^\prime(t),\boldsymbol{u}(t)\rangle_{\boldsymbol{V}^*\times\boldsymbol{V}},\;\;
\frac{d}{dt}\|\boldsymbol{\sigma}(t)\|^2_{\boldsymbol{L}_2}=2\langle\boldsymbol{\sigma}^\prime(t),\boldsymbol{\sigma}(t)\rangle_{\boldsymbol{L}_2^*\times\boldsymbol{L}_2},
\end{equation}
where 
$$
\boldsymbol{u}=\boldsymbol{v}_1-\boldsymbol{v}_2,\;\;\boldsymbol{\sigma}=\boldsymbol{\tau}_1-\boldsymbol{\tau}_2.
$$

Keeping in mind that $(\boldsymbol{v}_1,\boldsymbol{\tau}_1)$ and $(\boldsymbol{v}_2,\boldsymbol{\tau}_2)$ are weak solutions  to problem \eqref{equ-1}--\eqref{equ-6}, we obtain, after a simple calculation,
\begin{gather}
\text{Re}\langle\boldsymbol{u}^\prime,\boldsymbol{\varphi}\rangle-\text{Re}\sum\limits_{i=1}^{2}\Bigl(v_{1i}\boldsymbol{u},\frac{\partial{\boldsymbol{\varphi}}}{\partial x_i}\Bigr)
-\text{Re}\sum\limits_{i=1}^{2}\Bigl(u_i\boldsymbol{v}_2,\frac{\partial{\boldsymbol{\varphi}}}{\partial x_i}\Bigr)\nonumber
\\
+\bigl(\boldsymbol{\sigma},\boldsymbol{\mathcal E}(\boldsymbol{\varphi})\bigr)
+(1-a)\bigl(\nabla\boldsymbol{u},\nabla\boldsymbol{\varphi}\bigr)=0,\label{dp-2}
\end{gather}
\begin{equation}\label{dp-3}
\text{We}\langle\boldsymbol{\sigma}^\prime,\boldsymbol{\psi}\rangle+(\boldsymbol{\sigma},\boldsymbol{\psi})
=2a\bigl(\boldsymbol{\mathcal E}(\boldsymbol{u}),\boldsymbol{\psi}\bigr).
\end{equation}

Setting $\boldsymbol{\varphi}=\boldsymbol{u}(t)$ in \eqref{dp-2} and $\boldsymbol{\psi}=(2a)^{-1}\boldsymbol{\sigma}(t)$ in \eqref{dp-3}, we then add the results.
Taking into account \eqref{en-7}, we obtain
\begin{gather}
\frac{\text{Re}}{2}\frac{d}{dt}\|\boldsymbol{u}(t)\|^2_{\boldsymbol{L}_2}-\text{Re}\sum\limits_{i=1}^{2}\Bigl(u_i(t)\boldsymbol{v}_2(t),\frac{\partial{\boldsymbol{u}(t)}}{\partial x_i}\Bigr)\nonumber
\\
+\|\boldsymbol{u}(t)\|^2_{\boldsymbol{V}}
+\frac{\text{We}}{4a}\frac{d}{dt}\|\boldsymbol{\sigma}(t)\|^2_{\boldsymbol{L}_2}+\frac{1}{2a}\|\boldsymbol{\sigma}(t)\|^2_{\boldsymbol{L}_2}
=0.\label{en-9}
\end{gather}
Using an integration by parts, the H\"older inequality and the Ladyzhenskaya inequality, we get
\begin{gather}
\text{Re}\left|\sum\limits_{i=1}^{2}\Bigl(u_i(t)\boldsymbol{v}_2(t),\frac{\partial{\boldsymbol{u}(t)}}{\partial x_i}\Bigr)\right|
\leq C_3 \|\boldsymbol{u}(t)\|_{\boldsymbol{L}_4}^2\|\boldsymbol{v}_2(t)\|_{\boldsymbol{V}}\nonumber
\\
\leq C_4 \|\boldsymbol{u}(t)\|_{\boldsymbol{L}_2}\|\nabla\boldsymbol{u}(t)\|_{\boldsymbol{L}_2}\|\boldsymbol{v}_2(t)\|_{\boldsymbol{V}}
\leq C_5 \|\boldsymbol{u}(t)\|_{\boldsymbol{L}_2}\|\boldsymbol{u}(t)\|_{\boldsymbol{V}}\|\boldsymbol{v}_2(t)\|_{\boldsymbol{V}}\nonumber
\\
\label{en-10}
\leq C_6\|\boldsymbol{u}(t)\|_{\boldsymbol{L}_2}^2\|\boldsymbol{v}_2(t)\|_{\boldsymbol{L}_2}^2+\|\boldsymbol{u}(t)\|_{\boldsymbol{V}}^{2}.
\end{gather}
With \eqref{en-10} we deduce from \eqref{en-9} the estimate
$$
\frac{d}{dt}\left(\|\boldsymbol{u}(t)\|^2_{\boldsymbol{L}_2}+\|\boldsymbol{\sigma}(t)\|^2_{\boldsymbol{L}_2}\right)
\leq C_7\xi(t) \left(\|\boldsymbol{u}(t)\|^2_{\boldsymbol{L}_2}+\|\boldsymbol{\sigma}(t)\|^2_{\boldsymbol{L}_2}\right).
$$
where $\xi(t) =\|\boldsymbol{v}_2(t)\|_{\boldsymbol{V}}^2$,

 Note that  $\boldsymbol{u}(0)=\boldsymbol{0}$, $\boldsymbol{\sigma}(0)=\boldsymbol{0}$, and $\xi\in L_1(0,T)$. In this case, the Gr\"onwall--Bellman inequality gives
 $\boldsymbol{u}\equiv\boldsymbol{0}$ and $\boldsymbol{\sigma}\equiv\boldsymbol{0}$. Thus, we have proved that  problem \eqref{equ-1}--\eqref{equ-6} has a unique weak solution. 

By the above proof, we see that this solution belongs to $\boldsymbol{C}([0,T]; \boldsymbol{H}(\Omega))\times\boldsymbol{C}([0,T]; \boldsymbol{L}_2(\Omega))$. 
Moreover, we can derive the energy equation \eqref{eneq} in the standard way. This completes the proof of Theorem \ref{mainthr1}.

{\bf Received: July 29, 2016}

\end{document}